\documentclass[11pt]{amsart}

\usepackage{amssymb,amsmath}
\usepackage{verbatim}

\usepackage[pctex32]{graphics}

\theoremstyle{plain}
\newtheorem{Pro}{Proposition}[section]
\newtheorem{T}{Theorem}[section]
\newtheorem{Lemma}{Lemma}[section]
\newtheorem{C}{Corollary}[section]

\theoremstyle{definition}

\newtheorem{Example}{Example}[section]

\numberwithin{equation}{section}

\newcommand{\ko}{\hspace{2mm}}

\def \cal { \mathcal}

\def\qed{\hfill$\square$\smallskip}

\def\proof{\noindent{\sc Proof.}}

\begin{document}

\title[The first returning speed and the last exit speed  ]
{The first returning speed and the last exit speed of a type of Markov chain  }

\author{Huizeng Zhang*}
\thanks{*Department of mathematics, Hangzhou Normal University}
\thanks{*Research supported  by NSFC 11001070.}

\author{Minzhi Zhao** }
\thanks{**Department of Mathematics, Zhejiang University}
\thanks{**Research supported in part by NSFC 10601047  and  by Zhejiang Provincial Natural Science Foundation of China grant J20091364.}
\thanks{**Corresponding author.}

\author{Lei Wang***}
\thanks{ ***Department of mathematics, Zhejiang University}


\thanks{\it{2000 AMS Classification}: \rm{60J10}}
\begin{abstract}
Let $\{X_n\}$ be a Markov chain  with  transition probability $p_{ij}=a_{j-(i-1)^+},\forall i,j\ge 0$, where
$a_j=0$ provided $j<0$, $a_0>0$, $a_0+a_1<1$ and $\sum_{n=0}^\infty a_n=1$.  Let $\mu=\sum_{n=1}^\infty na_n$.
It's known that $\{X_n\}$ is positive recurrent when $\mu<1$; is null recurrent when $\mu=1$; and is transient
when $\mu>1$.  In this paper, we  shall discuss the first returning speed  and the last exit speed  more
precisely by means of $\{a_n\}$.
\end{abstract}
\maketitle


\section{Introduction}
Let $X=\{X_n:n=0,1,2,\cdots\}$ be a Markov chain with state space $\{0,1,\cdots\}$, transition matrix
$$P=\left(
      \begin{array}{ccccc}
       a_0& a_1 & a_2 & a_3 & \cdots \\
        a_0 & a_1 & a_2 & a_3 & \cdots \\
        0 & a_0 & a_1 & a_2 &\cdots \\
        0 & 0 & a_0 & a_1 & \cdots \\
        \vdots & \vdots & \vdots & \vdots & \vdots \\
      \end{array}
    \right).$$
That is,  $p_{ij}=a_{j-(i-1)^+}$, where $a_k=0$ for all $k<0$.
Suppose that $0<a_0<1$ ,  $a_0+a_1<1$ and $\sum\limits_{n=0}^\infty a_n=1$.  This Markov chain is known as model  repair shop, see \cite{KT} and \cite{B}.
The process $X$ is an irreducible Markov chain.
Let $\mu=\sum\limits_{n=1}^\infty na_n$.  Then $X$ is positive recurrent if $\mu<1$;   is null recurrent if $\mu=1$; and is transient  if $\mu>1$.

Suppose that $P(X_0=0)=1$,
 $\tau=\tau_0=\inf\{n\ge 1:X_n=0\}$ and $L=L_0=\sup\{n\ge 0:X_n=0\}$.  When $\mu=1$, $X$ is null recurrent. That is,
$P(\tau<\infty)=1$ and  $E(\tau)=\infty$. But for what $0<\alpha<1$ that can make $E(\tau^\alpha)<\infty$?
Similarly, when $\mu<1$, $X$ is positive recurrent. That is, $P(\tau<\infty)=1$ and $E(\tau)<\infty$.
But for what $1<\alpha <\infty$ that can make $E(\tau^\alpha)<\infty$? In a word, when $\mu\le 1$, the first returning
time to $0$ is finite almost surely, but how can we describe the returning speed in more detail?
When $\mu>1$, $X$ is transient, that is $P(L<\infty)=1$. On another words, the last exit time from $0$ is finite almost surely.
How can we characterize  the exit speed more precisely? For example,
 for what $0<\alpha <\infty$ that can make $E(L^\alpha)<\infty$?
In this paper, we shall answer these questions.

This paper is organized as follows. In \S2,
 the first returning speed and the last exit speed for any Markov
chain are discussed. Then in the following, we focus on the
 repair shop model Markov chain.
 In \S3,  the first returning speed  is discussed  in the case that $\mu=1$.
In \S4, the last exit speed is discussed in the case that $\mu>1$.
The first returning speed  in the case that  $\mu<1$ is discussed in \S5.

Throughout this paper, let $\mathbb N$  be the set of all positive integers.
For any sequence $\{w(n)\}$, let $\Delta w(n)=w(n)-w(n-1)$ and $\Delta^2 w(n)= \Delta w(n)-\Delta w(n-1)$.
Let $\cal W$  be the class of all nonnegative increasing sequence $(w(0),w(1),w(2),\cdots)$ such that
$ \Delta w(n)$ is decreasing and $\lim_{n\to\infty} \Delta w(n)=0$.
For any two nonnegative sequences $\{a_n\}$ and $\{b_n\}$, we use $a_n\thicksim b_n$ to denote the case that
there is $0<c_1<c_2<\infty$ and $M\in \mathbb N$ such that
$c_1a_n\le b_n\le c_2a_n$ for all $n\ge M$.
For any nonnegative  functions $f(x)$ and $g(x)$,
 $f(x)\thicksim g(x)$ as $x\to 0$ means that
 there is $0<c_1<c_2<\infty$ and $0<x_1<1$ such that
$c_1f(x)\le g(x)\le c_2 f(x)$ for all $0<x<x_1$.
Use $f(x)=o(g(x))$ as $x\to 0$ to denote the case that
$\lim_{x\downarrow 0}\frac {f(x)}{g(x)}=0$.
If $\lim_{x\downarrow 0}\frac {f(x)}{g(x)}=\infty$, then we write $f(x)=\Theta (g(x))$ as $x\to 0$.
Similarly $f(x)\thicksim g(x)$ as $x\to 1$ means  that
there is $0<c_1<c_2<\infty$ and $0<x_1<1$ such that
$c_1f(x)\le g(x)\le c_2 f(x)$ for all $x_1<x<1$.


\section{General results for Markov chains}
Suppose that $\{X_n\}$ is a  Markov chain on a probability space $(\Omega,\cal F, P)$, with $n$-step transition probabilities $p^{(n)}_{xy}$. Suppose that $x$ is a state and $P(X_0=x)=1$. Let $\tau$ and $L$ be the first returning time to $x$ and  the last exit time from $x$ respectively. That is
\begin{eqnarray*}
\tau&=&\inf\{n\ge 1:X_n=x\},\\
L&=&\sup\{n\ge 0: X_n=x\},
\end{eqnarray*}
where $\inf\emptyset=\infty$ and $\sup\emptyset=0$.
For any $t>0$, let
\begin{eqnarray*}
F(t)&=&E(t^{\tau};\tau<\infty\},\\
U(t)&=&\sum_{n=0}^\infty p^{(n)}_{xx}t^n,\\
L(t)&=&E(t^{L};L<\infty).
\end{eqnarray*}
It is known that
\begin{eqnarray}
U(t)=\left\{
         \begin{array}{ll}
           \frac 1{1-F(t)}, & \hbox{if $F(t)<1$;} \\
           \infty, & \hbox{otherwise.}
         \end{array}
       \right.  \label{eq2.1}
\end{eqnarray}
Let $p=F(1)=P(\tau<\infty), q=1-p$. Then for any nonnegative integer $n$,
$$P(L=n)=P(X_n=x,X_{n+1}\not=x,X_{n+2}\not=x,\cdots)=qP(X_n=x).$$
Hence we have
\begin{eqnarray}
 L(t)=qU(t).\label{e2.2}
\end{eqnarray}

The state $x$ is said to be recurrent if $p=1$, otherwise $x$ is said to be transient.
If $x$ is recurrent, then $P(L=\infty)=1$, otherwise $P(L<\infty)=1$.
That is, if $x$ is transient, then with probability one  $\{X_n\}$ will leave  from $x$  sooner or later.
Now a natural question has arisen: \textbf{how can we describe the last exit speed from $x$ more precisely? }
In \cite{SW}, the moments of last exit times for L\'evy processes are studied.
In \cite{ZZ}, the weighted moments of  last exit times for Markov processes are studied.
In this section, we shall use $F(t)$ and $U(t)$ to discuss the last exit speed for Markov chains.

On another hand, if $x$ is recurrent, then   though with probability one $\{X_n\}$ will return to $x$, but the first return speed may
be quite different.  To deal with this, recurrent states are divided into two classes: positive recurrent and null recurrent.
The sate $x$ is positive recurrent if and only if $E(\tau)<\infty$.
In \cite{Z}, we discussed the first returning speed of null recurrent Markov chains.
There, we proved that if $x$ is a null recurrent sate, then for any $w\in\cal W$,
\begin{eqnarray*}
E(w(\tau))<\infty \text{ if and only if } \sum_{n=1}^\infty -n \triangle^2 w(n+1)[F(1)-F(1-\frac 1n)]<\infty.
\end{eqnarray*}
But there are many positive recurrent Markov chains, and their returning speed
may be quite different too.  So
another question has arisen: \textbf{how can we characterize the first returning speed of positive recurrent states more precisely?}
In this section, we shall also discuss this question.

Let $R_0$ be the decay parameter of $x$, that is
$$R_0=\sup\{t\ge 0:F(t)\le 1\}.$$
Then by equation \ref{eq2.1},
$R_0=\sup\{t\ge 0: U(t)<\infty\}$.
If $F(R_0)=1$, equivalently, $U(R_0)=\infty$, then we say
$\{X_n\}$ is $R_0$-recurrent.
If $F'(R_0)<\infty$, that is, $E(R^\tau_0\tau;\tau<\infty)<\infty$, then $x$ is said to be $R_0$-positive recurrent.
In this case, $\{X_n\}$ has some nice properties, see for example \cite{J} and \cite{K}.
For discuss conveniently, we introduce another parameter
$$R_1=\sup\{t\ge 0: F(t)<\infty\}.$$
Obviously, $R_1\ge R_0\ge 1$. If $X$ is recurrent, then $R_0=1$. If $X$ is $R_0$-transient or $R_0$-null recurrent, then
$R_1=R_0$. The main results of this section are the following two  theorems. Theorem \ref{t2.1} is followed by lemma \ref{l2.2}
directly.

\begin{T}\label{t2.1}
 Suppose that $\{X_n\}$ is $R_0$-recurrent, $R_1<\infty$ and $F(R_1)<\infty$. Then we have the following properties.

(1) For any $k\in\mathbb N$,
$ E(R^{\tau}_1\tau^k;\tau<\infty)<\infty \text{\, if and only if\,} F^{(k)}(R_1)<\infty.$

(2) If there is an integer $k\ge 0$ such that $F^{(k)}(R_1)<\infty$ but  $F^{(k+1)}(R_1)=\infty$, then for any $w\in \cal W$,
$ E(R^{\tau}_1\tau^kw(\tau);\tau<\infty)<\infty$ if and only if
$\sum_{n=1}^\infty -n \Delta^2 w(n+1)[F^{(k)}(R_1)-F^{(k)}(R_1-\frac 1n)]<\infty$.
\end{T}

\begin{T}\label{t2.2}
Suppose that $\{X_n\}$ is $R_0$-transient. Then we have the following properties.

(1)That  $R_1=R_0$, $E(R^\tau_0;\tau<\infty)<1$ and $E(R^L_0)<\infty$.

(2) For any $k\in\mathbb N$, the following four conditions are equivalent to each other.

\hspace{1cm}(a)$E(R^{\tau}_0\tau^k;\tau<\infty)<\infty.$

\hspace{1cm}(b)  $F^{(k)}(R_0)<\infty$.

\hspace{1cm}(c) $U^{(k)}(R_0)<\infty$.

\hspace{1cm}(d) $E(R^{L}_0L^k)<\infty$.

(3) If there is an integer $k\ge 0$ such that $F^{(k)}(R_0)<\infty$ but  $F^{(k+1)}(R_0)=\infty$, then for any $w\in\cal W$,
the following four conditions are equivalent to each other.

(i) $E(R^{\tau}_0\tau^kw(\tau);\tau<\infty)<\infty$.

(ii)$E(R_0^{L}L^kw(L))<\infty$.

(iii)$\sum\limits_{n=1}^\infty -n \Delta^2 w(n+1)[F^{(k)}(R_0)-F^{(k)}(R_0-\frac 1n)]<\infty$.

(iv)$\sum\limits_{n=1}^\infty   -n \Delta^2 w(n+1)[U^{(k)}(R_0)-U^{(k)}(R_0-\frac 1n)]<\infty.$
\end{T}
\proof We need only to show  (2) and (3).

(2)
It suffices  to show that for any $k\in\mathbb N$,
$F^{(k)}(R_0)<\infty$ if and only if  $U^{(k)}(R_0)<\infty$.
The equation  $U(t)=1+U(t)F(t)$  implies that
\begin{eqnarray}
U^{(k)}(t)=\sum_{i=0}^k \binom {k}{i} U^{(i)}(t)F^{(k-i)}(t), \forall 0\le t\le R_0.\label{eq2.3}
\end{eqnarray}
It follows that
\begin{eqnarray}
U^{(k)}(R_0)(1-F(R_0))=U(R_0)F^{(k)}(R_0)+\sum_{i=1}^{k-1} \binom {k}{i} U^{(i)}(R_0)F^{(k-i)}(R_0).\label{eq2.4}
\end{eqnarray}
Particularly, $U'(R_0)(1-F(R_0))=U(R_0)F'(R_0)$. Since
$F(R_0)<1$ and $U(R_0)<\infty$, $U'(R_0)<\infty$ if and only if $F'(R_0)<\infty$.
Suppose that $k\ge 2$ and  we have showed that for any $n<k$,
$F^{(n)}(R_0)<\infty$ if and only if  $U^{(n)}(R_0)<\infty$.
By equation \ref{eq2.4}, $U^{(k)}(R_0)<\infty$ implies that $F^{(k)}(R_0)<\infty$.
Conversely, if $F^{(k)}(R_0)<\infty$, then
$F^{(n)}(R_0)<\infty$ for all $n<k$, and hence
$U^{(n)}(R_0)<\infty$ for all $n<k$. Since
$F(R_0)<1$ and $U(R_0)<\infty$, by (2.4), $U^{(k)}(R_0)<\infty$.
Thus our result holds by induction.

(3) By (2) and by Lemma \ref{l2.2}, it suffices  to show that for any  integer $k\ge 0$,
\begin{eqnarray}
U^{(k)}(R_0)-U^{(k)}(R_0-x)\thicksim F^{(k)}(R_0)-F^{(k)}(R_0-x), \text{\,as\,} x\to 0.\label{eq2.5}
\end{eqnarray}
If $k=0$, then
$$U(R_0)-U(R_0-x)=\frac {F(R_0)-F(R_0-x)}{(1-F(R_0))(1-F(R_0-x))}\thicksim F(R_0)-F(R_0-x), \text{\,as\,} x\to 0.$$
Now suppose that $k\ge 1$.
By equation \ref{eq2.4},
\begin{eqnarray*}
&&U^{(k)}(R_0)(1-F(R_0))-U^{(k)}(R_0-x)(1-F(R_0-x))\\
=&&[U(R_0)F^{(k)}(R_0)-U(R_0-x)F^{(k)}(R_0-x)]\\
 &&+\sum_{i=1}^{k-1} \binom {k}{i}
[U^{(i)}(R_0)F^{(k-i)}(R_0)-U^{(i)}(R_0-x)F^{(k-i)}(R_0-x)].
\end{eqnarray*}
Since $F'(R_0)<\infty$, $F(R_0)-F(R_0-x)\thicksim x$ as $x\to 0$.
Since $F^{(k+1)}(R_0)=\infty$, $U^{(k+1)}(R_0)=\infty$ and hence
$x=o(U^{(k)}(R_0)-U^{(k)}(R_0-x))$ as $x\to 0$. In addition $U^{(k)}(R_0)<\infty$.
Thus
\begin{eqnarray*}
&&U^{(k)}(R_0)(1-F(R_0))-U^{(k)}(R_0-x)(1-F(R_0-x))\\
=&&[U^{(k)}(R_0)-U^{(k)}(R_0-x)][1-F(R_0)]+U^{(k)}(R_0-x)[F(R_0-x)-F(R_0)]\\
\thicksim &&U^{(k)}(R_0)-U^{(k)}(R_0-x),  \text{\,as\,} x\to 0.
\end{eqnarray*}
Similarly, as $x\to 0$,
\begin{eqnarray*}
&&U(R_0)F^{(k)}(R_0)-U(R_0-x)F^{(k)}(R_0-x)\thicksim F^{(k)}(R_0)-F^{(k)}(R_0-x),\\
&&\sum_{i=1}^{k-1} \binom {k}{i}
[U^{(i)}(R_0)F^{(k-i)}(R_0)-U^{(i)}(R_0-x)F^{(k-i)}(R_0-x)]\thicksim x.
\end{eqnarray*}
The fact that $F^{(k+1)}(R_0)=\infty$ yields that $x=o(F^{(k)}(R_0)-F^{(k)}(R_0-x))$ as $x\to 0$.
Therefore \begin{eqnarray}
U^{(k)}(R_0)-U^{(k)}(R_0-x)\thicksim F^{(k)}(R_0)-F^{(k)}(R_0-x), \text{\,as\,} x\to 0.
\end{eqnarray}
\qed

Now we shall give  two lemmas that are used in  the proof of the theorems above.
Suppose that $\{a_0,a_1,a_2,\cdots \}$ is  a  nonnegative sequence.
For any positive integer $N$, let
$b_m^{(N)}=\sum_{n=mN}^{(m+1)N-1}a_n$ and $A_N=\sum_{n=0}^N a_n$.
Let  $A(t)=\sum_{n=0}^\infty a_n
t^n,  t\ge 0$.  Let $h(n)$ be an increasing positive
 sequence such that
$$0<\liminf_{n\to\infty}h(n)^n\le \limsup_{n\to\infty}h(n)^n<1.$$
The idea of Lemma \ref{l2.1} comes from Theorem 1 of \cite{IC}.

\begin{Lemma}\label{l2.1}
If there is $N_0$ and $x>0$ such that for all $N\ge
N_0$,
$$b_m^{(N)}\le xb^{(N)}_0, \forall m.$$
Then  $ A_n\sim A(h(n))$.
\end{Lemma}
\proof Obviously, $0<h(n)<1$ for all $n$. So we need only to
consider the case that $\sum_{n=0}^\infty a_n=\infty$. For any
$N$,
\begin{eqnarray*}
A(h(N))&=&\sum_{n=0}^\infty a_n h(N)^n\ge \sum_{n=0}^N
a_nh(N)^n\ge
h(N)^N\sum_{n=0}^N a_n\\
&=&A_Nh(N)^N.
\end{eqnarray*}
On  another hand, for any $N\ge N_0$,
\begin{eqnarray*}
A(h(N))&=&\sum_{n=0}^\infty a_n h(N)^n=\sum_{m=0}^\infty
\sum_{n=mN}^{(m+1)N-1}a_nh(N)^n\\
&\le & \sum_{m=0}^\infty h(N)^{mN}\sum_{n=mN}^{(m+1)N-1} a_n\le
\sum_{m=0}^\infty h(N)^{mN}(x\sum_{n=0}^N a_n)\\
&=&\frac {xA_N}{1-h(N)^N}.
\end{eqnarray*}
Hence $A_n\sim A(h(n))$.\qed

By Lemma \ref{l2.1},  if $c$ is a positive constant and  $\{a_n\}$ is a nonnegative decreasing sequence, then
 $A_n\sim A(1-\frac 1{cn})$.

Now suppose that $\sum\limits_{n=0}^\infty a_n<\infty$. Then for any positive integer $k$ and $0\le t<1$,
$$A^{(k)}(t)=\sum_{n=k}^\infty n(n-1)\cdots (n-k+1)a_nt^{n-k}.$$
Let
$$A^{(k)}(1)=\lim_{t\uparrow 1}A^{(k)}(t)=\sum_{n=k}^\infty  a_n n(n-1)\cdots
(n-k+1). $$ Then $\sum\limits_{n=0}^\infty a_n n^k<\infty$ if and only if
$A^{(k)}(1)<\infty$. Suppose that $\sum\limits_{n=0}^\infty {a_n}n^k<\infty$ and
$\sum\limits_{n=0}^\infty a_n n^{k+1}=\infty$. Let $w\in\cal W$. Now we show consider the finiteness
of $\sum\limits_{n=0}^\infty a_n n^k w(n)$.

\begin{Lemma}\label{l2.2}
Suppose that $c$ is a positive constant. Then $\sum\limits_{n=0}^\infty a_n n^k w(n)<\infty$ if and only if $\sum\limits_{n=1}^\infty -n \Delta^2
 w(n+1)[A^{(k)}(1)-A^{(k)}(1-\frac 1{cn})]<\infty$.
\end{Lemma}
\proof We may assume that $w(0)=0$ w.l.o.g. Let
$b_m=\sum_{n=m}^\infty a_nn^k$. Then $b_m\le b_0<\infty$ and
\begin{eqnarray*}
\sum_{n=0}^\infty  a_n n^k w(n)&=&\sum_{n=1}^\infty a_n n^k\sum_{m=1}^n \Delta w
(m)=\sum_{m=1}^\infty \Delta w(m)\sum_{n=m}^\infty a_n n^k\\
&=& \sum_{m=1}^\infty \Delta w(m)b_m=\sum_{m=1}^\infty
b_m\sum_{n=m}^\infty [\Delta w(n)-\Delta w(n+1)]\\
&=&\sum_{n=1}^\infty [\Delta w(n)-\Delta w(n+1)]\sum_{m=1}^n b_m.
\end{eqnarray*}

Since $b_m$ is decreasing, by Lemma \ref{l2.1},  $\sum_{m=1}^n b_m\sim B(1-\frac 1{cn})$,
where $B(t)=\sum_{n=1}^\infty b_n t^{n-1}$. On another hand,
\begin{eqnarray*}
\frac {A^{(k)}(1)-A^{(k)}(t)}{1-t}&=&\sum_{n=k+1}^\infty a_n
n(n-1)\cdots (n-k+1)(1+t+\cdots t^{n-k-1})\\
&=&\sum_{i=0}^\infty t^i\sum_{n=i+k+1}^\infty  a_n n(n-1)\cdots (n-k+1)\\
&\sim& \sum_{i=0}^\infty b_{i+k+1}t^i=\frac {B(t)}{t^k}[1-\frac{\sum_{j=1}^k b_j t^{j-1}}{B(t)}]
\text{\,as} \ \ t\to 1.
\end{eqnarray*}
Immediately,
$$B(1)=\sum_{m=1}^\infty b_m=\sum_{m=1}^\infty \sum_{n=m}^\infty a_nn^k=\sum_{n=1}^\infty a_nn^{k+1}=\infty.$$
Hence $\frac {A^{(k)}(1)-A^{(k)}(t)}{1-t}\thicksim B(t)$, as $t\to 1$. Therefore
 $\sum_{m=1}^n b_m\sim n[A^{(k)}(1)-A^{(k)}(1-\frac
1{cn})]$ which yields our result.\qed

By Lemma \ref{l2.2}, if there is $R>0$ and nonnegative integer $k$ such that  $\sum\limits_{n=0}^\infty {a_n}R^n n^k<\infty$ and
$\sum\limits_{n=0}^\infty a_n R^n n^{k+1}=\infty$, then for any  $w\in\cal W$,
 $\sum\limits_{n=0}^\infty a_n R^n n^k w(n)<\infty$ if and only if $\sum\limits_{n=1}^\infty -n \Delta^2
 w(n+1)[A^{(k)}(R)-A^{(k)}(R-\frac 1{n})]<\infty$.


\section{The first returning speed for null recurrent repair shop model}
In the following,
suppose that  $X=\{X_n:n=0,1,2,\cdots\}$ is  a repair shop model Markov chain.
That is,  $X$ is a Markov chain with state space $\{0,1,\cdots\}$,transition
  probability $p_{ij}=a_{j-(i-1)^+}$, where $a_k=0$ for all $k<0$,
 $0<a_0<1$ ,  $a_0+a_1<1$ and $\sum\limits_{n=0}^\infty a_n=1$.
Let $\mu=\sum\limits_{n=1}^\infty na_n$.
Suppose that $P(X_0=0)=1$ and $\tau=\inf\{n\ge 1:X_n=0\}$. For any $t>0$, let
 $G(t)=\sum\limits_{n=0}^\infty a_n t^n$ and $F(t)=E(t^{\tau};\tau<\infty)$.     Then by \cite{KT}, for any $t>0$,
\begin{eqnarray}
F(t)=tG(F(t)).\label{eq3.1}
\end{eqnarray}

 In this section we shall consider the case that $\mu=1$,  that is,  $X$ is null recurrent.

Firstly, we shall characterize the asymptotic  behavior of $1-F(t)$ as $t\to 1$ by means   of $G(t)$.
Let $\mathcal H$ be the collection of all nonnegative strictly increasing continuous functions  on $[0,1]$.
For any $f\in\cal H$, the  inverse function of $f$ exists, namely $f^{-1}$, defined on $[f(0),f(1)]$.
 For convenience, we may extend the domain of $f^{-1}$ to $[f(0),\infty)$ by defining
$f^{-1}(x)=1$ for all $x\ge f(1)$. Then $f^{-1}$ is strictly increasing on $[f(0),f(1)]$, and increasing on $[f(0),\infty)$.

For any $h\in [0,1]$, let
$$\phi(h)=1-G'(1-h), \psi(h)=\int_0^h\phi(x)\,dx=G(1-h)-(1-h).$$
Since $\phi(h)>0$ for all $0<h<1$, $\psi(h)$ is a strictly increasing continuous function on $[0,1]$, that is
$\psi\in\cal H$.
Since $\psi(0)=0$ and $\psi(1)=a_0$, 
 $\psi^{-1}$ is strictly increasing on $[0,a_0]$,  increasing on $[0,\infty)$,  $\psi^{-1}(0)=0$ and $\psi^{-1}(x)=1$
 for all $x\ge a_0$.

Let $\mathcal H'$ be the collection of all $f\in\mathcal H$ that satisfying the following two properties:

(c1) For all $x,y\ge 0$ with $x+y\le 1$, $ f(x+y)\ge f(x)+f(y)$.

(c2) For any positive integer $n$ and $x\in[0,\frac 1n]$,
$f(nx)\le n^2 f(x).$

For example, $f(x)=x^\alpha$ is in $\mathcal H'$ if and only if $1\le \alpha\le 2$.
We have the following lemmas.
\begin{Lemma}\label{l3.1}
(1)That $\psi\in\cal H'$.

(2)For any $f\in\cal H'$,$f^{-1}$ is a sub-additive function on $[0,\infty)$, that is for any
    $x,y\ge 0$, $$f^{-1}(x+y)\leq f^{-1}(x)+f^{-1}(y).$$

(3) For any $f\in\cal H'$, any positive integer $n$ and any $0\le y\le f(\frac 1n)$,
 $$f^{-1}(n^2 y)\ge  n f^{-1}(y).$$

\end{Lemma}
\proof (1) For any $h\in (0,1)$, $\phi'(h)=G''(1-h)>0$, so $\phi$ is strictly increasing on $[0,1]$.
Thus for any $h_1,h_2\ge 0$ with $h_1+h_2\le 1$,
\begin{eqnarray*}
 \psi(h_1+h_2)&=&\int^{h_1}_{0}\phi(x)dx+\int^{h_1+h_2}_{h_1}\phi(x)dx\\
&\geq& \int^{h_1}_{0}\phi(x)dx+\int^{h_2}_{0}\phi(x)dx=\psi(h_1)+\psi(h_2).
\end{eqnarray*}
Similarly, since $\phi''(h)=-G^{(3)}(1-h)\le 0$ on $(0,1)$, $\phi'(h)$ is decreasing on $[0,1]$ and hence $\phi(h)$ is a
sub-additive function on $[0,1]$. Thus for any $x\in [0,\frac 1n]$,
$$\psi(nx)=\int_0^x n\phi(nh)dh\le n^2\int_0^x \phi(h)dh=n^2\psi(x).$$
Consequently,  $\psi\in\cal H'$.

(2) Suppose that $f\in\cal H'$. Then $f(0+0)\ge f(0)+f(0)$ implies that $f(0)=0$.
Hence the domain of $f^{-1}$ is $[0,\infty)$, $f^{-1}(0)=0$, $f^{-1}(x)=1$ for all $x\ge f(1)$.
For any
 $x,y\ge 0$,
if $f^{-1}(x)+f^{-1}(y)\ge 1$, then
$f^{-1}(x+y)\leq f^{-1}(x)+f^{-1}(y)$ holds.

If $f^{-1}(x)+f^{-1}(y)<1$, then
$f[f^{-1}(x)+f^{-1}(y)]\ge f[f^{-1}(x)]+f[f^{-1}(y)]=x+y$.
It follows that $f^{-1}(x)+f^{-1}(y)\ge f^{-1}(x+y)$.

Therefore (2)holds.

(3) Let $x=f^{-1}(y)$. Then $0\le x\le \frac 1n$ and  $f(nx)\le n^2 f(x)=n^2 y$. Hence
$f^{-1}(n^2 y)\ge n f^{-1}(y)$ holds.
\qed

\begin{Lemma}\label{l3.2}
 For any $0<t< 1$, $1-F(t)=\psi^{-1}(\frac {(1-t)F(t)}{t})$.
\end{Lemma}
\proof
For any  $0< t< 1$, by equation (3.1),
$$\psi(1-F(t))=G(F(t))-F(t)=
\frac {(1-t)F(t)}{t}$$ which implies that
$1-F(t)=\psi^{-1}(\frac {(1-t)F(t)}{t})$.\qed

\begin{T}\label{t3.1}
As $t\rightarrow 1$,$1-F(t)\thicksim \psi^{-1}(1-t)$.
\end{T}
\proof
Since $\lim_{t\to 1}\frac {F(t)}{t}=1$, there is $0<t_0<1$, such that for any $t_0<t<1$,
$$\frac 12 \le \frac {F(t)}{t}\le 2.$$
Thus for any $t_0<t<1$,
$$\psi^{-1}[\frac 12 (1-t)]\le \psi^{-1}[\frac {F(t)}{t}(1-t)]\le \psi^{-1}[2(1-t)].$$
By Lemma \ref{l3.1} and Lemma \ref{l3.2}, for any $t_0<t<1$,
$$\frac 12 \psi^{-1}(1-t)\le 1-F(t)\le 2\psi^{-1}(1-t).$$
Consequently, as $t\rightarrow 1$,$1-F(t)\thicksim \psi^{-1}(1-t)$.
\qed

Particularly, $1-F(1-\frac 1n)\thicksim \psi^{-1}(\frac 1n)$. By Theorem \ref{t2.1}, we have the following corollary.

\begin{C}\label{c3.1}
(1)For any $w\in\cal W$,
 $E[w(\tau)]<\infty$ if and only if $\sum\limits_{n=1}^\infty -n\Delta^2 w(n+1)\psi^{-1}(\frac 1n)<\infty$.

(2) For any $0<\alpha<1$, $E(\tau^\alpha)<\infty$ if and only if $\sum\limits_{n=1}^\infty n^{\alpha-1}\psi^{-1}(\frac 1n)<\infty$.
\end{C}
Corollary \ref{c3.1} tells us that the returning speed of $0$ can be characterized by $\psi$.
The quicker that $\psi^{-1}(x)$ converges to $0$  as $x$
tends to $0$, the quicker that $\{X_n\}$ returns to $0$. Since $\psi^{-1}$ is dependent on $G(x)$ completely, we
can judge the finiteness of $E[w(\tau)]$ by means of $G(x)$ directly. Next, we shall discuss the returning speed in more detail.

\begin{Lemma}\label{l3.3}
 Suppose that $\psi_1,\psi_2\in\cal H'$ and $n$ is a positive integer.

(1)  If $\psi_1(x)\le n\psi_2(x)$ as $x\to 0$,  then
$\psi_2^{-1}(y)\le n \psi_1^{-1}(y)$ as $y\to 0$.

(2) If $\psi_1(x)\ge n^2\psi_2(x)$ as $x\to 0$,  then
$\psi_2^{-1}(y)\ge n \psi_1^{-1}(y)$ as $y\to 0$.
\end{Lemma}
\proof   (1) There is $0<x_0<1$ such that  $\psi_1(x)\le n\psi_2(x)$ for all $0<x<x_0$.
Let $y_0=\psi_1(x_0)$. Then $y_0>\psi_1(0)=0$. For any $0<y<y_0$, $0<\psi_1^{-1}(y)<x_0$, and hence
$$\psi_1(\psi_1^{-1}(y))\le n\psi_2(\psi_1^{-1}(y)).$$  That is,  $y\le n\psi_2(\psi_1^{-1}(y))$.
Consequently,
$$\psi^{-1}_2(y)\le \psi^{-1}_2[n\psi_2(\psi_1^{-1}(y))]\le n \psi^{-1}_2[\psi_2(\psi_1^{-1}(y))]=n\psi_1^{-1}(y).$$

(2) There is $0<x_0<\frac 1n$ such that $\psi_1(x)\ge n^2\psi_2(x)$  for all $0<x<x_0$.
Thus by Lemma \ref{l3.1}, for any $0<x<x_0$,
$$\psi^{-1}_2[\psi_1(x)]\ge \psi^{-1}_2[n^2\psi_2(x)]\ge n \psi^{-1}_2[\psi_2(x)]=nx.$$
Letting $y=\psi_1(x)$, we get that $\psi^{-1}_2(y)\ge n\psi^{-1}_1(y)$ holds for all
$0<y<\psi_1(x_0)$.
\qed

By Lemma \ref{l3.3} and Theorem \ref{t3.1}, we immediately get the following proposition.

\begin{Pro} \label{p3.1}
Suppose that $\psi_1\in\cal  H'$.

(1) If as $x\to 0$, $\psi(x)=o(\psi_1(x))$, then as $t\to 1$,  $1-F(t)=\Theta(\psi_1^{-1}(1-t))$

(2) If as $x\to 0$, $\psi(x)=\Theta(\psi_1(x))$, then as $t\to 1$, $1-F(t)=o(\psi_1^{-1}(1-t))$.

(3) If as $x\to 0$, $\psi_1(x)\thicksim\psi(x)$, then as $t\to 1$, $1-F(t)\thicksim \psi_1^{-1}(1-t)$.

\end{Pro}

\begin{C}\label{c3.2}
(1)If $G''(1)<\infty$, then as $t\to 1$,  $1-F(t)\thicksim\sqrt{1-t}$,  and $E(\tau^\alpha)<\infty$ if and only if $\alpha<\frac 12$.

(2) If $G''(1)=\infty$, then as $t\to 1$, $1-F(t)=o(\sqrt{1-t})$, and $E(\tau^\alpha)<\infty$ for all $\alpha<\frac 12$.

(3)If there is $0<\beta<1$ such that  $1-G'(x)\thicksim(1-x)^\beta$ as $x\to 1$, then as $t\to 1$, $1-F(t)\thicksim(1-t)^{\frac 1{1+\beta}}$,
and $E(\tau^\alpha)<\infty$ if and only if $\alpha<\frac 1{\beta+1}$.
\end{C}
\proof (1) Since $\lim_{x\uparrow 1}\frac {1-G'(x)}{1-x}=G''(1)<\infty$,
as $h\to 0$, $\phi(h)=1-G'(1-h)\thicksim h$.
Thus as $x\to 0$, $$\psi(x)=\int_0^x \phi(h)dh\thicksim \int_0^x hdh \thicksim x^2.$$
It is easy to verify that the function $\psi_1(x):=x^2$ is in $\cal H'$. By Proposition \ref{p3.1},
as $t\to 1$, $1-F(t)\thicksim \psi^{-1}_1(1-t)=\sqrt{1-t}$.
By Theorem \ref{t2.1}, $E(\tau^\alpha)<\infty$ if and only if $\alpha<\frac 12$.

(2)Since $G''(1)=\infty$, as $h\to 0$, $\phi(h)=1-G'(1-h)=\Theta(h)$. It implies that
as $x\to 0$, $$\psi(x)=\int_0^x \phi(h)dh=\Theta(x^2).$$
By Proposition  \ref{p3.1} and Theorem  \ref{t2.1}, our result holds.

(3) As $h\to 0$, $\phi(h)=1-G'(1-h)\thicksim h^\beta$. Thus as $x\to 0$, $$\psi(x)=\int_0^x \phi(h)dh\thicksim \int_0^x h^\beta dh \thicksim x^{1+\beta}.$$
It is easy to verify that the function $\psi_2(x)=x^{1+\beta}$ is in $\cal H'$. By Proposition \ref{p3.1},
as $t\to 1$, $1-F(t)\thicksim \psi^{-1}_2(1-t)=(1-t)^{\frac 1{1+\beta}}$.
By Theorem \ref{t2.1}, $E(\tau^\alpha)<\infty$ if and only if $\alpha<\frac 1{1+\beta}$.\qed

Proposition \ref{p3.1} and Corollary \ref{c3.2}  tell us that the quicker that $G(x)$ tends to $1$ as $x$ tends to $1$, the more slowly  that
$\{X_n\}$ returns to $0$. In the case that  $G''(1)<\infty$,  $\{X_n\}$ returns to $0$ the most slowly. In this case $E(\tau^\alpha)<\infty$
if and only if $\alpha<\frac 12$.

\begin{Example}\label{ex3.1}
 Suppose that $a_n=\frac 1{2^{n+1}},n\ge 0$. Then $\mu=1$ and $\{X_n\}$ is null recurrent.
It is easy to check that  $G(t)=\frac 1{2-t}$.
Thus $F(t)=tG(F(t))=\frac {t}{2-F(t)}$, which implies that $F(t)=1-\sqrt{1-t}$.
Particularly, $1-F(1-\frac 1n)=\frac 1{\sqrt n}$.
Therefore $E(\tau^\alpha<\infty)$ if and only if $\alpha<\frac 12$.

We may deduce this result by Corollary \ref{c3.2} directly. In fact
$G'(t)=\frac 1{(2-t)^2}$ and $G''(t)=\frac 2{(2-t)^3}$. Hence $G''(1)=2<\infty$. So by Corollary \ref{c3.2},
$1-F(t)\thicksim \sqrt{1-t}$, and $E(\tau^\alpha)<\infty$ if and and only if $\alpha<\frac 12$.
\end{Example}

\begin{Example}\label{ex3.2}
 Suppose that $a_0=\frac 23$, $a_1=0$, $a_2=\frac 14$, and for $n\ge 3$,
$$a_n=\frac {2}{3}{(-1)^n} \binom {\frac 32}{n}=\frac {1\times 3 \times 5\times\cdots\times(2n-5)}{2^{n-1}n!}>0.$$
Then $G(x)=\sum\limits_{i=0}^\infty  a_i x^i=x+\frac 23 (1-x)^{\frac 32}$ and $G'(x)=1-\sqrt{1-x}$.
Thus $\sum\limits_{n=0}^\infty a_n=G(1)=1$ and $\mu=G'(1)=1$.
So $\{X_n\}$ is a null recurrent Markov chain.
By $F(t)=t(G(F(t))$, we get $\frac {(1-t)F(t)}{t}=\frac 23(1-F(t))^\frac 32$ which implies that
$1-F(t)\thicksim (1-t)^\frac 23$.
So $E(\tau^\alpha)<\infty$
if and only if  $\alpha<\frac 23$.

We may also deduce this result by Corollary \ref{c3.2} directly.
In fact, since   $1-G'(x)\thicksim (1-x)^{\frac 12}$ as $x\to 1$, by Corollary \ref{c3.2}, $1-F(t)\thicksim (1-t)^{\frac 23}$
as $t\to 1$.
\end{Example}
\section{The last exit speed for transient repair shop model}

In this section, suppose that $\mu>1$. Then $\{X_n\}$ is transient. We shall discuss the last exit speed from $0$.
Let $L=\sup\{n\ge 0:X_n=0\}$, $R_0$ be the decay parameter of $\{X_n\}$.
If $F(R_0)=1$, then $F(R_0)=R_0G(F(R_0))$ implies that $R_0=1$.
It contradicts that $F(1)<1$. So $F(R_0)<1$ and $X$ is $R_0$-transient.
In addition, $F(t)=\infty$ provided $t>R_0$. For any $x\ge 0$ with $G'(x)<\infty$,
let $$\xi(x)=G(x)-xG'(x).$$
Then $\xi(0)=a_0>0$, $\lim_{t\uparrow 1}\xi(1)=1-\mu<0$, and  $\xi'(x)=-xG''(x)<0$ for all $x>0$.
So there exits unique $0<x_0<1$ such that $\xi(x_0)=0$, that is,  $G(x_0)=x_0G'(x_0)$.

\begin{T}\label{t4.1}
(1) That $R_0=\frac {x_0}{G(x_0)}>1$,  $F(R_0)=x_0<1$ and $E(R_0^L)<\infty$.

(2)That  $E(R_0^\tau\tau,\tau<\infty)=E(R_0^L L)=\infty$.

(3) For any $w\in\cal W$, $E(R_0^\tau w(\tau);\tau<\infty)<\infty$ if and only if
$E(R_0^L w(L))<\infty$, and if and only if $\sum_{n=1}^\infty -n \Delta^2 w(n+1)\psi_{x_0}^{-1} (\frac 1n)<\infty$,
where $\psi_{x_0}=\frac {G(x_0(1-h))}{G(x_0)}-(1-h)\in\cal H'$ and $ \psi_{x_0}^{-1} $ is the inverse function of
 $\psi_{x_0}$.
\end{T}
\proof \ko (1) At first, we shall introduce a transformation.  Suppose that $x>0$ and $G(x)<\infty$.
Define $P^{(x)}(X_0=0):=1$.
For any $n\ge 1$ and any nonnegative integers $i_0=0,i_1,\cdots,i_n$, define
\begin{eqnarray*}
P^{(x)}(X_0=0,X_1=i_1,\cdots,X_n=i_n):&=&\frac {x^{\sum_{k=1}^n i_k-(i_{k-1}-1)^{+}}}{G(x)^n}P(X_0=0,X_1=i_1,\cdots,X_n=i_n)\\
&=& \prod_{k=1}^n \frac {x^{i_k-(i_{k-1}-1)^+}a_{i_k-(i_{k-1}-1)^+}}{G(x)}.
\end{eqnarray*}
Let $\cal F_n=\sigma\{X_0,X_1,\cdots,X_n\}$ and $\cal F=\sigma\{X_0,X_1,\cdots\}$.
It is easy to verify that $P^{(x)}$ is a pre-probability measure on $\bigcup_n \cal F_n$.
Therefore $P^{(x)}$ has a unique extension to $\cal F$ which is  denoted  still by $P^{(x)}$. Then
$P^{(x)}$ is a probability measure on $\cal F$. Under $P^{(x)}$, $\{X_n\}$ is a homogenous Markov chain
with  transition probability $$p^{(x)}_{ij}=\frac {x^{j-(i-1)^+}a_{j-(i-1)^+}}{G(x)}.$$
Let
\begin{equation}\label{eq4.1}
a^{(x)}_j= \frac {a_j x^{j}}{G(x)},\forall j\ge 0.
\end{equation}
Then $\sum_{j=0}^\infty a^{(x)}_j=\sum_{j=0}^\infty \frac {a_j x^{j}}{G(x)}=1$,
$a^{(x)}_0>0, a^{(x)}_0+a^{(x)}_1=1-\sum_{j=2}^\infty a^{(x)}_j<1$.
Thus under $P^{(x)}$, $\{X_n\}$ is  a repair shop model with $\{a^{(x)}_j\}$.

Next we shall calculate the value $F(t)$.
For any $t>0$,
let
\begin{eqnarray}
G^{(x)}(t):&=&\sum_{j=0}^\infty a^{(x)}_j t^j=\frac {G(xt)}{G(x)},\label{eq4.2}\\
F^{(x)}(t):&=&E^{(x)}(t^{\tau};\tau<\infty)=\sum_{n=1}^\infty P^{(x)}(\tau=n)t^n.\label{eq4.3}
\end{eqnarray}
For any $n\ge 1$ and any nonnegative integers $i_0=0,i_1\not=0,i_2\not=0,\cdots,i_{n-1}\not=0,i_n=0$,
\begin{eqnarray*}
P^{(x)}(X_0=0,X_1=i_1,\cdots,X_n=i_n)= \frac {x^{n-1}}{G(x)^n}P(X_0=0,X_1=i_1,\cdots,X_n=i_n).
\end{eqnarray*}
It follows that
\begin{eqnarray}
P^{(x)}(\tau=n)=\frac {x^{n-1}}{G(x)^n}P(\tau=n),\forall n \label{eq4.4}
\end{eqnarray}
 which implies   that
$$F^{(x)}(t)=\frac {F(\frac {xt}{G(x)})}{x}.$$
Particularly, $F(\frac x{G(x)})=x F^{(x)}(1)\le x$.

Let $\mu^{(x)}=\sum\limits_{j=1}^\infty j a^{(x)}_j$. Then
\begin{equation}
\mu^{(x)}=\frac {xG'(x)}{G(x)}.\label{eq4.5}
\end{equation}
 If $\mu^{(x)}\le 1$, then
$\{X_n\}$ is recurrent under $P^{(x)}$, that is, $F^{(x)}(1)=1$. Thus
if $ \frac {xG'(x)}{G(x)}\le 1$, then $F(\frac x{G(x)})=x$. For any $x\ge 0$ with $G(x)<\infty$,
let $$\eta(x)=\frac x{G(x)}.$$
 Then $\eta'(x)=\frac {\xi(x)}{G(x)^2}$ for all $x$ satisfying $G'(x)<\infty$.
Thus $\eta'(x)>0$ on $(0,x_0)$, and $\eta'(x)<0$ when $x>x_0$. It leads to that
$\eta(x)$ is strictly increasing on $[0,x_0)$, and then strictly   decreasing when $x>x_0$. The maximum of
$\eta(x)$ is $\eta(x_0)$. For any $y\ge 0$ with $F(y)<\infty$, $y=\frac {F(y)}{G(F(y))}=\eta(F(y))\in [0,\eta(x_0)]$.
That is, $F(y)<\infty $ if and only if
$y\le \frac {x_0}{G(x_0)}$.
For any $y\in[0,\eta(x_0)]$, there is unique $x\in [0,x_0]$ such that
$y=\eta(x)$. Since $\eta'(x)\ge 0$,   $ \frac {xG'(x)}{G(x)}\le 1$ and $F(y)=x$, that is
$F(y)$ is the minimum solution of  $\eta(x)=y$.

By the discussion above, we get that $R_0= \eta(x_0)=\frac {x_0}{G(x_0)}$ and $F(R_0)=x_0<1$.
Since $x_0\not =1$, $R_0=\eta(x_0)>\eta(1)=1$. By Theorem \ref{t2.2}, $E(R_0^L)<\infty$.

(2) By equation \ref{eq4.5}, $\mu^{(x_0)}=\frac {x_0G'(x_0)}{G(x_0)}=1$. It follows that $\{X_n\}$ is null
recurrent under $P^{(x_0)}$.  That is,  $P^{(x_0)}(\tau<\infty)=1$ and $E^{(x_0)}(\tau)=\infty$. By Corollary \ref{c3.1},
for any $w\in\cal W$, $E^{(x_0)}(w(\tau))<\infty$ if and only if
$\sum_{n=1}^\infty -n \Delta^2 w(n+1)\psi_{x_0}^{-1}(\frac 1n)<\infty$, where
$\psi_{x_0}(h)=G^{(x_0)}(1-h)-(1-h)$. By equation \ref{eq4.2}, $\psi_{x_0}(h)=\frac {G(x_0(1-h))}{G(x_0)}-(1-h)$.
By equation \ref{eq4.4},
$E^{(x_0)}(\tau)=\frac {E(R_0^\tau \tau;\tau<\infty)}{x_0}$ and
$E^{(x_0)}(w(\tau))=\frac {E(R_0^\tau w(\tau);\tau<\infty)}{x_0}$ for all $w\in \cal W$.
Combination with Theorem \ref{t2.2},
 (2) and (3) hold.\qed

\begin{Example}\label{ex4.1}
 Suppose that  $a_n=q^np, n\ge 0$, where $0<p<\frac 12$ and $q=1-p$. Then $G(t)=\frac p{1-tq}$ and $\mu=\frac qp>1$. So $\{X_n\}$ is transient.
 By equation $F(t)=tG(F(t))$, we get that $F(t)=\frac {1-\sqrt{1-4pqt}}{2q}$.  Thus   $R_0=\frac 1{4pq}$,
$F(R_0)=\frac 1{2q}<1$ and hence $\{X_n\}$ is $R_0$-transient. Further $F'(R_0)=\infty$ and
$F(R_0)-F(R_0-\frac 1n)=\frac 1{2q}\sqrt{\frac 1{nR_0}}$. Hence by Theorem \ref{t2.2}, for any $\alpha>0$,
$E( R_0^{L}L^\alpha)<\infty$ if and only if $\alpha<\frac 12$.

We may use Theorem \ref{t4.1} to deduce the results  above. Clearly,  $G'(t)=\frac {pq}{(1-tq)^2}$. Thus
the equation $G(x)=xG'(x)$ has unique solution $x_0=\frac 1{2q}$. Hence $R_0=\frac {x_0}{G(x_0)}=\frac 1{4pq}$.
In addition
$a_n^{(x_0)}=\frac {a_nx_0^n}{G(x_0)}=\frac 1{2^{n+1}}, n\ge 0$.
That is  under $P^{(x_0)}$, $\{X_n\}$ is the same process as in Example \ref{ex3.1}.
Therefore for any $\alpha>0$,
$E( R_0^{L}L^\alpha)<\infty$ if and only if $\alpha<\frac 12$.
\end{Example}

\section{The first returning speed for positive recurrent repair shop model}

In this section,  suppose that $\mu<1$. Then $\{X_n\}$ is positive recurrent. We shall discuss the returning time of $0$.
Let $$R_1=\sup\{t:F(t)<\infty\},$$
 $$R=\sup\{t:G(t)<\infty\}.$$
Let $D=\{x\ge 0: G(x)<\infty\}$ and $D'=\{x\ge 0:G'(x)<\infty\}$.
If $R=\infty$, then $D=D'=[0,\infty)$.
If $R<\infty$ and $G'(R)<\infty$, then $G(R)<\infty$ and  $D=D'=[0,R]$.
If $R<\infty$, $G(R)<\infty$ but $G'(R)=\infty$, then $D=[0,R]$ and $D'=[0,R)$.
If $R<\infty$ and $G(R)=\infty$, then $D=D'=[0,R)$.

\begin{T}\label{t5.1}
(1) If $G'(R)<\infty$ and $G(R)>RG'(R)$, then $R_1=\frac R{G(R)}$ and $F(R_1)=R$.

(2) If $G'(R)=\infty$ or $G(R)\le RG'(R)$, then there is unique $x_0>1$ such that
$G(x_0)=x_0G(x_0)$. In this case $R_1=\frac {x_0}{G(x_0)}>1$ and $F(R_1)=x_0>1$.

(3) That $R_1=1$ if and only if $R=1$.
\end{T}
\proof
For any $x\in D$, let
 $$\eta(x)=\frac {x}{G(x)}.$$
For any $x\in D'$, let
$$\xi(x)=G(x)-xG'(x).$$
Then $\eta'(x)=\frac {\xi(x)}{G^2(x)}$ for all $x\in D'$.
Since $\xi(0)=a_0>0$ and $\xi'(x)=-xG''(x)<0$, $\xi(x)=0$ has at most one solution.
If $\xi(x)=0$ has no solution, then $\xi>0$ on $D'$ and $\eta$ is an  increasing function on $D$.
If $\xi(x)=0$ has a unique solution namely $x_0$, then $\xi>0$ on $[0,x_0)$ and $\xi<0$ on $(x_0,\infty)\cap D'$,
which follows that $\eta$ is strictly increasing on $[0,x_0]$ and then strictly decreasing when $x>x_0$.
By the proof of Theorem \ref{t4.1},
$F(y)<\infty$ if and only if $y=\eta(x)$ for some $x$ with $G(x)<\infty$, and
in this case,
$$F(y)=\{x:\eta(x)=y,G(x)\ge x G'(x)\}=\{x:\eta(x)=y, \eta'(x)\ge 0\}=min\{x:\eta(x)=y\}.$$

If $\xi(x)=0$  has a unique  solution $x_0$,
then the range of $\eta(x)$ is $[0,\eta(x_0)]$. So $F(t)<\infty$ if and only if $t\le \eta(x_0)$.
That is $R_1=\frac {x_0}{G(x_0)}$. Also $F(R_1)=x_0$. Since $\xi(1)=1-\mu>0$, $x_0>1$ and $R_1=\eta(x_0)>\eta(1)=1$.

If $\xi(x)=0$ has no  solutions, then   the range of $\eta(x)$ is $[0,\eta(R)]$ if $G(R)<\infty$,
or $[0,\eta(R))$ if $G(R)=\infty$. We shall show that $G'(R)<1$. If $R=1$, then $G'(1)=\mu<1$ holds.
Now suppose that $R>1$. Note that $\xi(x)>0$ for all $x<R$ and $\eta(x)$ is strictly increasing on $[0,R]$.
For any $\frac {1+R}2<x<R$,
$$G'(x)<\frac {G(x)}{x}=\frac 1{\eta(x)}\le \frac 1{\eta(\frac {1+R}{2})}.$$
Thus
$$G'(R)=\lim_{x\uparrow R} G'(x)\le \frac 1{\eta(\frac {1+R}{2})}<\frac 1{\eta(1)}=1.$$
Therefore $G'(R)<1$ and $G(R)<\infty$.
Consequently, $R_1=\eta(R)=\frac R{G(R)}$ and $F(R_1)=R$.
If $R=1$, then $R_1=\eta(1)=1$. If $R>1$, then $R_1=\eta(R)>\eta(1)=1$.
\qed

Similar as the proof of the statements (2) and (3) in Theorem \ref{t4.1}, we may prove the following theorem.

\begin{T}\label{t5.2}
Suppose that  there is unique $x_0>1$ such that
$G(x_0)=x_0G'(x_0)$. Then
$E(R_1^\tau\tau)=\infty$, and for
 any $w\in\cal W$, $E(R_1^\tau w(\tau))<\infty$ if and only if
 $\sum_{n=1}^\infty -n \Delta^2 w(n+1)\psi_{x_0}^{-1} (\frac 1n)<\infty$,
where $\psi_{x_0}=\frac {G(x_0(1-h))}{G(x_0)}-(1-h)\in\cal H'$ and $ \psi_{x_0}^{-1} $ is the inverse function of
 $\psi_{x_0}$.
\end{T}

We shall discuss the first returning speed of $\{X_n\}$ further in the case that $G(x)=xG'(x)$ has no positive solutions.
\begin{Lemma}\label{l5.1}
 For any integer $n\ge 2$, $E(\tau^n)<\infty$ if and only if $G^{(n)}(1)<\infty$.
\end{Lemma}
\proof  By $F(t)=tG(F(t))$, we get that
\begin{eqnarray}
F'(t)&=&G(F(t))+tG'(F(t))F'(t),\label{eq5.1}\\
F''(t)&=&2G'(F(t))F'(t)+tG''(F(t))[F'(t)]^2+tG'[F(t)]F''(t).\label{eq5.2}
\end{eqnarray}
We shall use induction to show that
\begin{eqnarray}
F^{(n)}(t)=H_{n-1}(t)+tG^{(n)}[F(t)][F'(t)]^n+tG'[F(t)]F^{(n)}(t), \forall n\ge 2,\label{eq5.3}
\end{eqnarray}
where $H_{n-1}(t)$ is a polynomial function of $F'(t),...,F^{(n-1)}(t),G'(F(t)),...,G^{(n-1)}(F(t))$
with positive coefficients.

By  equation \ref{eq5.2}, equation \ref{eq5.3} holds for $n=2$.
Suppose that equation  \ref{eq5.3} holds for some  $n=k\ge 2$, that is
$$F^{(k)}(t)=H_{k-1}(t)+tG^{(k)}[F(t)][F'(t)]^k+tG'[F(t)]F^{(k)}(t).$$
Then $$F^{(k+1)}(t)=H_{k}(t)+tG^{(k+1)}[F(t)][F'(t)]^{k+1}+tG'[F(t)]F^{(k+1)}(t),$$ where
\begin{eqnarray*}
H_{k}(t)&=&H'_{k-1}(t)+G^{(k)}[F(t)][F'(t)]^k+ktG^{(k)}[F(t)][F'(t)]^{k-1}F''(t)\\
&&+G'[F(t)]F^{(k)}(t)+tG''[F(t)]F'(t)F^{(k)}(t).
\end{eqnarray*}
Obviously, $H_{k}(t)$ is a polynomial function of  $F'(t),...,F^{(k)}(t),G'(F(t)),...,G^{(k)}(F(t))$ with positive coefficients.
Thus equation \ref{eq5.3} holds.

It follows that
\begin{eqnarray}
F^{(n)}(t)=\frac{H_{n-1}(t)+tG^{(n)}[F(t)][F'(t)]^n}{1-tG'[F(t)]}, \forall n\ge 2.\label{eq5.4}
\end{eqnarray}
Particularly, by equation \ref{eq5.1} and \ref{eq5.4},
\begin{eqnarray}
F'(1)&=&\frac 1{1-\mu},\label{eq5.5}\\
F^{(n)}(1)&=&\frac {H_{n-1}(1)+G^{(n)}(1)(\frac 1{1-\mu})^n}{1-\mu}, \forall n\ge 2.\label{eq5.6}
\end{eqnarray}
Obviously,
$F''(1)<\infty$ if and only if $G''(1)<\infty$.
Suppose that  when $n\le k$, $F^{(n)}(1)<\infty$ if and only if $G^{(n)}(1)<\infty$.

If $F^{(k+1)}(1)<\infty$, then $F^{(n)}(1)<\infty$ for all $n\le k$, and hence $G^{(n)}(1)<\infty$ for all
$n\le k$, which yields that $H_{k}(1)<\infty$. Consequently $G^{(k+1)}(1)<\infty$.
Similarly, if $G^{(k+1)}(1)<\infty$, then $F^{(k+1)}(1)<\infty$.

Therefore by induction, for all $n\ge 2$, $F^{(n)}(1)<\infty$ if and only if
$G^{(n)}(1)<\infty$.  Now the fact that $E(\tau^n)<\infty$ if and only if $F^{(n)}(1)<\infty$ yields our result.
\qed

\begin{Lemma} \label{l5.2}
Suppose that there is $k\in\mathbb  N$ such that
$E(\tau^k)<\infty$ but $E(\tau^{k+1})=\infty$. Then for any $w\in\cal W$,
$E(\tau^k w(\tau))<\infty$ if and only if $\sum_{n=1}^\infty  n^kw(n)a_n<\infty$.
\end{Lemma}

\proof  By Lemma \ref{l5.1}, $G^{(k)}(1)<\infty$ but $G^{(k+1)}(1)=\infty$.
By Theorem \ref{t2.1}, for any $w\in\cal W$,
$E(\tau^k w(\tau))<\infty$ if and only if  $\sum_{n=1}^\infty -n\Delta^2w(n+1) [F^{(k)}(1)-F^{(k)}(1-\frac 1n)]<\infty$.
By Lemma \ref{l2.2}, for any $w\in\cal W$,
$\sum_{n=1}^\infty  n^kw(n)a_n<\infty$ if and only if $\sum_{n=1}^\infty -n\Delta^2w(n+1) [G^{(k)}(1)-G^{(k)}(1-\frac 1n)]<\infty$.

Thus it suffices to show that
\begin{eqnarray}
F^{(k)}(1)-F^{(k)}(1-t)\thicksim G^{(k)}(1)-G^{(k)}(1-t), \text{\,as\,} t\to 0.\label{eq5.7}
\end{eqnarray}
Firstly, we shall show that
\begin{eqnarray}
G^{(k)}(1)-G^{(k)}[F(1-t)]\thicksim G^{(k)}(1)-G^{(k)}(1-t), \text{\,as\,} t\to 0.\label{eq5.8}
\end{eqnarray}
For $0<t<1$,
 $F(1-t)\le 1-t$ and hence $G^{(k)}(1)-G^{(k)}[F(1-t)]\ge G^{(k)}(1)-G^{(k)}(1-t)$.
Since the function $G^{(k)}(x)$ is convex,
\begin{eqnarray*}
G^{(k)}(1)-G^{(k)}[F(1-t)]\le [1-F(1-t)]\frac {[G^{(k)}(1)-G^{(k)}(1-t)]}{t}\le F'(1) [G^{(k)}(1)-G^{(k)}(1-t)].
\end{eqnarray*}
Consequently, \ref{eq5.8}  holds.

If $k=1$, then
\begin{eqnarray*}
F'(1)-F'(1-t)=\frac {1-(1-t)G'[F(1-t)]-(1-\mu)G[F(1-t)]}  {(1-\mu)\{1-(1-t)G'[F(1-t)]\}}.
\end{eqnarray*}
Obviously,
\begin{eqnarray*}
&& 1-(1-t)G'[F(1-t)]-(1-\mu)G[F(1-t)]\\
&=&(1-\mu)\{G[F(1)]-G[F(1-t)]\}+\{\mu-G'[F(1-t)]\}+tG'[F(1-t)].
\end{eqnarray*}
Note that  $G'(1)<\infty$ and $F'(1)<\infty$, combination with \ref{eq5.8}, we get that  as $t\to 0$,
\begin{eqnarray*}
(1-\mu)\{G[F(1)]-G[F(1-t)]\}&\thicksim& t,\\
\mu-G'[F(1-t)]&\thicksim& G'(1)-G'(1-t)\\
tG'[F(1-t)]&\thicksim& t,
\end{eqnarray*}
and $$(1-\mu)\{1-(1-t)G'[F(1-t)]\}\thicksim (1-\mu)^2.$$
The fact that $G''(1)=\infty$ implies that $t=o(G'(1)-G'(1-t))$ as $t\to 0$.
Therefore as $t\to 0$,
 $F'(1)-F'(1-t)\thicksim G'(1)-G'(1-t)$.

If $k\ge 2$, then by equation \ref{eq5.3},
\begin{eqnarray*}
F^{(k)}(1)-F^{(k)}(1-t)
&=&\frac {H_{k-1}(1)-H_{k-1}(1-t)}{1-(1-t)G'[F(1-t)]}\\
&+&\frac {\{\mu-(1-t)G'[F(1-t)]\}[H_{k-1}(1)+G^{(k)}(1)(\frac 1{1-\mu})^k]} {(1-\mu)\{1-(1-t)G'[F(1-t)]\}}\\
&+&\frac {G^{(k)}(1)(\frac 1{1-\mu})^k-(1-t)G^{(k)}[F(1-t)][F'(1-t)]^k} {1-(1-t)G'[F(1-t)]}
\end{eqnarray*}
By \ref{eq5.8} and by the fact that $t\thicksim o(G^{(k)}(1)-G^{(k)}(1-t))$ as $t\to 0$, we get
\begin{eqnarray*}
&&G^{(k)}(1)(\frac 1{1-\mu})^k-(1-t)G^{(k)}[F(1-t)][F'(1-t)]^k\\
&&=(1-t)(\frac 1{1-\mu})^k\{G^{(k)}(1)-G^{(k)}[F(1-t)]\}+(1-t)G^{(k)}[F(1-t)]\{(\frac 1{1-\mu})^k-[F'(1-t)]^k\}\\
&&\,\,+tG^{(k)}(1)(\frac 1{1-\mu})^k\\
&&\thicksim G^{(k)}(1)-G^{(k)}(1-t), \text{\, as \,} t\to 0.
\end{eqnarray*}
In addition,  as $t\to 0$,
\begin{eqnarray*}
&&H_{k-1}(1)-H_{k-1}(1-t)\thicksim t,\\
&&\{\mu-(1-t)G'[F(1-t)]\}[H_{k-1}(1)+G^{(k)}(1)(\frac 1{1-\mu})^k]\thicksim t,\\
&&1-(1-t)G'[F(1-t)]\thicksim 1-\mu.
\end{eqnarray*}
Therefore $$F^{(k)}(1)-F^{(k)}(1-t)\thicksim G^{(k)}(1)-G^{(k)}(1-t)  \text{\,as\,} t\to 0,$$
and our result holds.
\qed

\begin{T}\label{t5.3}
 If $G'(R)<\infty$ and $G(R)>RG'(R)$, then we have the following properties.

(1)That  $E(R_1^\tau \tau)<\infty$ ,

(2)  For any $\alpha>1$, $E(R_1^\tau \tau^\alpha)<\infty$ if and only if $\sum_{n=1}^\infty  R^n n^\alpha a_n<\infty$.

(3)  If there is $k\in \mathbb N$ such that  $E(R_1^\tau \tau^k)<\infty$ but  $E(R_1^\tau \tau^{k+1})=\infty$,
then for any $w\in\cal W$, $E(R_1^\tau \tau^k w(\tau))<\infty$ if and only if $\sum_{n=1}^\infty R^n n^k  w(n)a_n<\infty$.
\end{T}
\proof Since $\mu^{(R)}=\frac {RG'(R)}{G(R)}<1$, $\{X_n\}$ is still positive recurrent under
$P^{(R)}$. Hence $E^{(R)}(\tau)<\infty$. By Lemma \ref{l5.1}, for any integer $k\ge 2$,
$E^{(R)}(\tau^k)<\infty$ if and only if $\sum_{n=1}^\infty a_n^{(R)}n^k<\infty$.  In addition, by Lemma \ref{l5.2}, if there is $k\in \mathbb N$ such that  $E^{(R)}(\tau^k)<\infty$ but  $E^{(R)}(\tau^{k+1})=\infty$, then for any $w\in\cal W$,
$E^{(R)}(\tau^k w(\tau))<\infty$ if and only if $\sum_{n=1}^\infty  n^kw(n)a^{(R)}_n<\infty$. On another hand,
$E^{(R)}(\tau^k)=\frac {E(R_1^\tau \tau^k)}{R}$, $E^{(R)}(\tau^kw(\tau))=\frac {E(R_1^\tau \tau^k w(\tau))}{R}$ and
$a_n^{(R)}=\frac {a_n R^n}{G(R)}$.
Thus our result holds. \qed

\begin{Example}\label{ex5.1}
Suppose that  $a_n=q^np, n\ge 0$, where $\frac 12<p<1$ and $q=1-p$. Then $G(t)=\frac p{1-tq}$ and $\mu=\frac qp<1$. So  $\{X_n\}$ is positive recurrent.  Similar as Example \ref{ex4.1},
  $F(t)=\frac {1-\sqrt{1-4pqt}}{2q}$.  Thus   $R_1=\frac 1{4pq}$, $F'(R_1)=\infty$ and  $F(R_1)-F(R_1-\frac 1n)=\frac 1{2q}\sqrt{\frac 1{R_1n}}$.
So by Theorem \ref{t2.1}, $E(R_1^\tau\tau^\alpha)<\infty$ if and only if $\alpha<\frac 12$.

We may use Theorem \ref{t5.1} and Theorem \ref{t5.2} to deduce our result.
Immediately,
$R=\sup\{t:G(t)<\infty\}=\frac 1q$ and $G'(R)=\infty$. Thus the
equation $G(x)=xG'(x)$ has unique solution, indeed, it is $x_0=\frac 1{2q}$. Thus
$R_1=\frac {x_0}{G(x_0)}=\frac 1{4pq}$. Also  under $P^{(x_0)}$, $\{X_n\}$ is the same process as in Example \ref{ex3.1}.
Therefore for any $\alpha>0$,
$E(R_1^\tau\tau^\alpha)<\infty$ if and only if $\alpha<\frac 12$.
\end{Example}

\begin{Example}\label{ex5.2}
suppose that  $\alpha>2$ and
 $a_k=(\frac 1{k+1})^\alpha-(\frac 1{k+2})^\alpha, k\ge 0$.
Then
$$\sum_{k=0}^\infty k a_k=(\frac 1 2)^\alpha+(\frac 13)^\alpha+\cdots\le \int_1^\infty (\frac 1x)^\alpha dx=\frac 1{\alpha-1}<1.$$
Thus $\{X_n\}$ is a positive recurrent. As $k\to\infty$,
$a_k\thicksim k^{-\alpha-1}$.
So $G(t)=\infty$ for all $t>1$ and $R=1$.   By Theorem \ref{t5.1},  $R_1=1$.   For any
$\beta>0$, $\sum_{k=1}^\infty a_kk^\beta<\infty$  if and only if
$\beta<\alpha$. So by Theorem \ref{t5.3},  $E(\tau^\beta)<\infty$ if and only if $\beta<\alpha$.
\end{Example}



\bigskip
\begin{thebibliography}{99}
\addcontentsline{toc}{chapter}{BIBLIOGRAPHY}


\bibitem{B} Pierre  Bremaud, \sc Markov Chains: Gibbs Fields, Monte Carlo Simulation, and Queues, \rm Springer-Verlag New York(1999).


\bibitem{IC} Dean Isaacson and Peter Colwell, \it Levels of Null
Persistency for Markov Chains, \rm J.Appl.Prob, 19,425-429(1982).


\bibitem{J} D. Vere-Jones, \it Geometric ergodicity in denumerable Markov Chains, \rm Quart. J. Math. Vol., 13, No. 2, 7-28(1962).

\bibitem{K} J. F. C. Kingman, \it The exponential decay of Markov transition probabilities, \rm Proc. London Mth. Soc., Vol. 13,
No.3, 337-358(1963).


\bibitem{KT} S. Karlin and H. M. Taylor, \sc A first course in stochastic processes, \rm Academic Press(1975).



\bibitem{SW}  Ken-iti Sato, and Toshiro WatanabeAnnales,\it   Moments of last exit times for L¨¦vy processes, \rm Ann. Inst. Henri Poincar\'e, Probab. Statist., Volume 40, Issue 2, 207-225(2004).

\bibitem{Z} M. Z.  Zhao, \it The first returning speed of null recurrent Markov chians,
\rm Chineses Ann.  Math. Ser. A, Vol.  27, No. 6, 761-770(2006).


\bibitem{ZZ} M. Z. Zhao and H. Z. Zhang, \it The weighted transience and recurrence of Markov processes, \rm Acta Mathematic
Sinica,English Seiries,  Vol. 23, No. 1, 111-126(2007).


\end{thebibliography}
\end{document}